\documentclass{amsart}
\usepackage{amscd,amssymb}

\newcounter{noindnum}
\newcommand{\noindstep}{\refstepcounter{noindnum}{\rm(}\alph{noindnum}\/{\rm)}}

\newcommand{\smcn}{{\rm;}}
\newcommand{\RightHamRed}{\,\mathop{\mathstrut_0\!\!\!\setminus\!\!\setminus}}

\renewcommand{\phi}{\varphi}
\renewcommand{\epsilon}{\varepsilon}
\renewcommand{\emptyset}{\varnothing}

\newcommand{\C}{\mathbb C}
\newcommand{\divisor}{\mathfrak D}
\newcommand{\ad}{\mathop{\mathrm{ad}}\nolimits}
\newcommand{\Ad}{\mathop{\mathrm{Ad}}\nolimits}

\newcommand{\E}{\mathcal{E}}
\newcommand{\F}{\mathcal{F}}
\newcommand{\fh}{\mathfrak h}
\newcommand{\fg}{\mathfrak g}
\newcommand{\fu}{\mathfrak u}
\newcommand{\fb}{\mathfrak b}
\newcommand{\fhr}{{\mathfrak h^r}}
\newcommand{\fgr}{{\mathfrak g^{rs}}}

\newcommand{\dual}{\mathbb I}
\newcommand{\rk}{\mathop{\mathrm{rk}}}
\newcommand{\Spec}{\mathop{\mathrm{Spec}}}
\newcommand{\Res}{\mathop{\mathrm{Res}}}
\newcommand{\Lie}{\mathop{\mathrm{Lie}}}

\newcommand{\conn}{\mathop{\mathbf{Conn}}\nolimits}
\newcommand{\Conn}{\mathcal{C}onn}
\newcommand{\fConn}{\mathcal{C}onn_{\mathrm{fr}}}
\newcommand{\overfConn}{\overline{\mathcal{C}onn}_{\mathrm{fr}}}
\newcommand{\CircX}{{\mathop{X}\limits^\circ}}
\newcommand{\OpenX}{{\dot X}}
\newcommand{\const}{\mathop{\mathrm{const}}}
\newcommand{\Hyper}{\mathbb H}

\theoremstyle{plain}
\newtheorem{theorem}{Theorem}
\newtheorem{proposition}{Proposition}
\newtheorem*{corollary}{Corollary}
\newtheorem{lemma}{Lemma}

\theoremstyle{definition}
\newtheorem*{convention}{Convention}
\newtheorem{definition}{Definition}

\theoremstyle{remark}
\newtheorem{remark}{Remark}

\title[Isostokes Deformations]{Algebraic and hamiltonian approaches\\
to isostokes deformations}

\author{Roman M.\ Fedorov}
\thanks{Partially supported by NSF grant DMS-0401164}
\address{The University of Chicago, Department of Mathematics,
5734~S. University~Ave, Chicago, Illinois 60615, USA}
\email{fedorov@math.uchicago.edu}

\begin{document}

\begin{abstract}
We study a generalization of the isomonodromic deformation to the case of
connections with irregular singularities. We call this generalization
Isostokes Deformation. A new deformation parameter arises: one can deform the
formal normal forms of connections at irregular points. We study this part of
the deformation, giving an algebraic description. Then we show how to use
loop groups and hypercohomology to write explicit hamiltonians. We work on an
arbitrary complete algebraic curve, the structure group is an arbitrary
semisimiple group.
\end{abstract}

\maketitle

\thispagestyle{empty}

\tableofcontents{} \vfill\eject

\section*{Introduction}

The isomonodromic deformation is a~classical subject pertaining to many areas
of mathematics (see~\cite{Babelon, Its} for example). In~\cite{Jimbo} it was
generalized to the case of connections with arbitrary order poles. In this
case one requires the monodromy data and the Stokes multipliers to remain
constant. Thus we suggest the term \emph{Isostokes Deformation}. In irregular
case a~new direction of deformation arises: one can deform the irregular
types of connections at irregular singular points. Thus one can deform the
curve, the divisor, and the irregular types. The deformation of the curve and
the divisor was further studied in~\cite{Krichever} and~\cite{FrenkelBenZvi}.
We study the deformation of the irregular types, our approach is close to
that of~\cite{FrenkelBenZvi}.

The deformation of the irregular types was also studied
in~\cite{BoalchThesis}. In that paper the algebraic curve is $\C P^1$, the
structure group is $GL(n)$. Our interest in this subject was evoked by the
papers~\cite{Boalch} and~\cite{FrenkelBenZvi}. In the former the deformation
for the divisor $2(0)+(\infty)$ is studied for an arbitrary complex reductive
group $G$. Its monodromy turns out to coincide with the action of
a~generalized braid group on the dual Poisson group $G^*$. The hamiltonian
approach in this case is obtained in~\cite{Harnad}. It turns out that in this
case the isostokes connection is the quasi-classical limit of the De
Concini--Milson--Toledano Laredo (or DMT) connection (see~\cite{Laredo,
Laredo2}). A conjecture of De Concini and Toledano Laredo says that the
monodromy of the DMT connection coincides with the action of the braid group
on a quantum group (this conjecture has been proved recently by Toledano
Laredo, see~\cite{Laredo3}). Thus the result of~\cite{Boalch} can be thought
as a geometrization of the quasi-classical limit of this conjecture.

We give an algebraic description of the isomonodromic deformation. Then we
give a~hamiltonian description with explicit hamiltonians. As by-products we
obtain a description of algebraic and Poisson structures on the moduli spaces
of connections.

The structure of the paper is as follows. In the first section we present
main definitions and results. We generalize the notion of analytic isostokes
deformation to arbitrary complete smooth complex algebraic curves and
arbitrary complex semisimple groups (the precise definition is given in \S2,
see Proposition~\ref{analyticisostokes}). To give an algebraic description of
the deformation we put a~structure of algebraic stack on the moduli space of
connections. The result of Proposition~\ref{modulistack} looks classical but
we could not find any reference. Similar constructions are discussed
in~\cite{Arinkin}.

Then we give an algebraic description (Theorem~\ref{algan}) of the isostokes
deformation. To obtain a~hamiltonian description we define the moduli stack
of connections with unipotent structures (see \S\ref{shamappr} and
Proposition~\ref{connu}). Finally, we give a~hamiltonian description
(Theorem~\ref{hamappr}).

In the rest of the paper we prove these theorems. In section~\ref{sheur} we
give an explanation of Theorem~\ref{hamappr} via loop groups --- this is how
the theorem was invented. This is an explanation in the spirit
of~\cite{FrenkelBenZvi} we were looking for. Another proof is given in
\S\ref{sehamappr}. In the last section we count the dimensions and explain
how to write explicit formulae.

The author wants to thank D.~Arinkin for invaluable numerous discussions from
which he learnt a~lot. The idea of Theorem~\ref{algan} belongs to him.
Without his help this paper would never be finished. I am grateful to my
advisor V.~Ginzburg and to D.~Ben--Zvi for many useful comments and
discussions. I am also thankful to V.~Baranovsky, E.~Frenkel, D.~Gaitsgory,
and V.~Vologodsky for their interest in my work.

\section{Main results}

\subsection{Bundles with connections}

\subsubsection{Non-resonant connections} Let us fix a~smooth
complete algebraic curve $X$ over $\C$, an effective divisor
$\divisor=\sum_1^ln_ix_i$ ($n_i>0$) on $X$, and a~connected semisimple group
$G$ over $\C$. We shall call $x_i$ an \emph{irregular point\/} if $n_i\ge2$.
We assume that there is at least one irregular point.

Fix an analytic coordinate $z_i$ at $x_i$ for $i=1,\ldots,l$. These
coordinates will be fixed throughout the paper. We could have avoided fixing
coordinates, then we would have to work with jets. However, it would make
things more complicated.

Consider a~pair $(E,\nabla)$, where $E$ is a~principal $G$-bundle on $X$ with
a~\emph{left\/} action of $G$, $\nabla$ is a~singular connection on $E$ such
that the polar divisor of this connection is bounded by $\divisor$. Choose
any trivialization of $E$ in the formal neighborhood of $x_i$. This
trivialization allows us to identify connections on the restriction of $E$ to
this neighbourhood with formal $\fg$-valued 1-forms, where $\fg$ is the Lie
algebra of $G$. Thus we can write:
\begin{equation}\label{leading}
\nabla=d+A_{n_i}\frac{dz_i}{z_i^{n_i}}+O\left(z_i^{1-n_i}\right),
\end{equation}
where $A_{n_i}\in\fg$. Denote by $\fgr$ the set of regular semisimple
elements of $\fg$. The connection $\nabla$ is called \emph{non-resonant\/} at
an irregular point $x_i$ if ${A_{n_i}\in\fgr}$; the connection $\nabla$ is
called \emph{non-resonant\/} if it is non-resonant at all the irregular
points.

The conjugacy class of $A_{n_i}$ does not depend on the choice
of trivialization of $E$ in the neighborhood of $x_i$. Thus the notion of
non-resonant connection does not depend on the choice of trivialization. We
shall call pairs $(E,\nabla)$ \emph{connections\/} for brevity.

We shall denote by $\overline{\Conn}$ the moduli space of pairs $(E,\nabla)$,
where $E$ is a~principal $G$-bundle, $\nabla$ is a~connection on $E$ with the
polar divisor bounded by $\divisor$. Let $\Conn$ be the subspace of
$\overline{\Conn}$ corresponding to the non-resonant connections. We shall
see below that $\overline{\Conn}$ has a~natural structure of an algebraic
stack, clearly, $\Conn$ is its open substack (see
Proposition~\ref{modulistack}).

\begin{remark}
Notice that in~\cite{Boalch} and elsewhere it is customary to write
a~connection as $d-A(z)$, since one thinks about a connection as about
a~differential equation. We always write connections in the form $d+A(z)$.
Notice also that in some papers on the subject (in particular,
in~\cite{Boalch}) $G$-bundles with \emph{right\/} actions are considered.
\end{remark}

\subsubsection{Compatible framings}\label{families} Let $x_i$ be an irregular
point. A \emph{framing\/} of $E$ at $x_i$ is a~choice of an element $s_i$ in
the fiber of $E$ over $x_i$.

We shall fix a~maximal torus $T$, a~maximal unipotent subgroup $U$ and
a~Borel subgroup $B$ in $G$ such that $T\subset B$, $U\subset B$. Let $\fh$,
$\fu$ and $\fb$ be the corresponding Cartan, maximal nilpotent, and Borel
subalgebras respectively.

After a framing at $x_i$ is chosen, the coefficient $A_{n_i}$ at the leading
term in~(\ref{leading}) is well-defined (not up to a conjugation). The
framing is called \emph{compatible\/} with~$\nabla$ if $A_{n_i}\in\fh$.
Compatible framings at $x_i$ for a~non-resonant connection $\nabla$ form an
$N(T)$-torsor, where $N(T)$ is the normalizer of $T$ in $G$. Of course, not
every resonant connection has a~compatible framing.

Denote by $\overfConn$ the moduli space of triples $(E,\nabla,s)$,
where $(E,\nabla)\in\overline{\Conn}$, $s=\{s_i\}$ is a~collection of
compatible framings at all the irregular points. There is a~natural forgetful
map $\overfConn\to\overline{\Conn}$, the preimage of $\Conn$ will be denoted
by $\fConn$. Clearly, $\fConn$ is an open subspace of $\overfConn$.

\begin{proposition}\label{modulistack}
The moduli space $\overline{\Conn}$ is an algebraic stack, while $\overfConn$
is an algebraic space.
\end{proposition}

We shall prove it in~\S\ref{secstack}.

\subsection{Formal normal forms of connections at irregular singular points}
Set $\fhr=\fh\cap\fgr$. Let $\nabla$ be a~non-resonant connection, $x_i$ be
an irregular point, $s_i$ be a~compatible framing at $x_i$. Then we can
choose a~trivialization of $E$ in the formal neighborhood of $x_i$ such that
$\nabla$ takes its formal normal form (see~\cite{Boalch}, Lemma~1):
\begin{equation}\label{fnf}
\nabla=d+h_{n_i}\frac{dz_i}{z_i^{n_i}}+
h_{n_i-1}\frac{dz_i}{z_i^{n_i-1}}+\ldots+h_1\frac{dz_i}{z_i},\qquad n_i\ge2,
\end{equation}
where $h_{n_i}\in\fhr$, $h_j\in\fh$ for $j=n_i-1,\ldots,1$. Notice that
without the framing this normal form would be defined  up to the diagonal
action of the Weyl group of $\fg$.

Taking the formal normal form at every irregular point $x_i$ we get a~map:
\[
\fConn\to(\fhr)^{l_{irr}} \times(\fh)^{\deg\divisor-l},
\]
where $l_{irr}$ is the number of irregular points, $l$ is the total number of
singular points of~$\divisor$, $\deg\divisor$ is the sum of multiplicities of
singular points.

Consider the map obtained from the previous map by forgetting the formal
residue~$h_1$ at every irregular point $x_i$ (in other words, this map
assigns to a~connection its \emph{irregular type}):
\begin{equation}\label{IT}
IT:\fConn\to(\fhr)^{l_{irr}}\times(\fh)^{\deg\divisor-l-l_{irr}}.
\end{equation}
This map will be of primary interest for us.

\subsection{Analytic isostokes deformations}
\begin{convention}
Let $\Delta$ be an algebraic scheme or a~smooth analytic manifold. By
a~\emph{non-resonant family of connections\/} over $\Delta$ we mean
a~triple $(E(t),\nabla(t),s(t))$, where $t\in\Delta$ is the deformation
parameter, $E(t)$ is a~principal $G$-bundle on $X\times\Delta$, $\nabla(t)$
is a~connection on $E(t)$ along $X$ with the polar divisor bounded by
$\divisor\times\Delta$ such that $\nabla(t)$ is non-resonant on the fiber
over any point of $\Delta$, $s(t)$ is a set of compatible framings at all the
irregular points. A~\emph{framing\/} for a~family $E(t)$ at $x_i\in X$ is
a~section of $E(t)|_{\{x_i\}\times\Delta}$.
\end{convention}

In~\S\ref{sanis} we shall construct a~natural connection on the map
(\ref{IT}) in the following sense. Given $(E,\nabla,s)\in\fConn$ and an
analytic map $f$ from a~polydisk $\Delta$ to the target space of $IT$ such
that
\begin{equation}\label{chp}
IT(E,\nabla,s)=f(0),
\end{equation}
we produce a~canonical way to extend $(E,\nabla,s)$ to a~non-resonant family
of connections over $\Delta$. This will be called \emph{Isostokes
Deformation}. Heuristically, we deform the connection in such a~way that both
monodromy data and Stokes data remain constant.

\begin{remark}
We think about the fibers of $IT$ as about generalized topological data
associated to connections. Thus, loosely speaking, the isostokes deformation
is the deformation of non-topological data (i.e.\ of the irregular type),
while preserving the topological data. More generally, one can consider
a~deformation preserving generalized topological data but changing irregular
types, the curve $X$, and the divisor $\divisor$. An approach based on loop
groups to the deformation of the curve and the divisor is given in
\cite{FrenkelBenZvi}.
\end{remark}

\subsection{Algebraic approach to isostokes deformations}\label{ssalgap}
Consider deformations that do not change a given connection but change
framings. It will be convenient for us to regard such deformations as
isostokes.

Let $\Delta\ni t_0$ be a smooth manifold, $v\in T_{t_0}\Delta$,
$(E(t),\nabla(t),s(t))$ be a non-resonant family, parameterized by $\Delta$
($T_{t_0}\Delta$ is the tangent space to $\Delta$ at $t_0$). To give an
algebraic description of the isostokes deformations we use a notion of a
family \emph{infinitesimally isostokes in the direction of $v$\/}.
Intuitively, it means that the restriction of this family to
$\dual=\Spec\C[\epsilon]/\epsilon^2$ is isostokes (where we view $v$ as a map
$\dual\to\Delta$) and we shall give a precise definition in~\S\ref{sinis}.

Let $X'$ be a~subset of $X$. We say that \emph{the restriction of
$(E(t),\nabla(t),s(t))$ onto~$X'$ is algebraically (analytically) constant\/}
if this restriction is algebraically (analytically) isomorphic to the
pullback of a~triple $(E_0,\nabla_0,s_0)$ along the projection
$X'\times\Delta\to X'$. In our applications $X'$ will be either open or
closed in $X$.

\begin{definition}\label{isostokesvectors}
Consider the open algebraic curve $\OpenX=X\setminus\divisor$ and let $v$ be
a tangent vector to $\fConn$ at $(E,\nabla,s)$. As a~tangent vector to the
moduli space it induces an \emph{algebraic\/} family of connections over
$\dual$. We call $v$ \emph{isostokes\/} if the restriction of this family to
$\OpenX\times\dual$ is algebraically constant.
\end{definition}

For any map $f$ we denote the corresponding tangent map by $f_*$.

\begin{theorem}\label{algan}
Consider an algebraic non-resonant family $(E(t),\nabla(t),s(t))$
parameterized by a~smooth variety $\Delta\ni t_0$, $v\in T_{t_0}\Delta$. Let
$f:\Delta\to\fConn$ be the induced map to the moduli space. The family is
infinitesimally isostokes at $t_0$ in the direction of $v$ iff $f_*v$ is an
isostokes vector.
\end{theorem}

\begin{remark}
Notice that since we work on the open curve, an analytically constant family
is not always algebraically constant. Moreover, roughly speaking, the theorem
above states that the Stokes data does not change in a family, whose
restriction to $\OpenX$ is algebraically constant. Of course, it is not true
for a family with the analytically constant restriction.
\end{remark}

\subsection{Hamiltonian approach}\label{shamappr}

Now we want to give a~hamiltonian description of the isostokes deformation.
Actually, $\fConn$ is a~Poisson space. It is not hard to see that connections
that are in the same symplectic leaf of $\fConn$ have the same irregular
type. Hence, our deformation is transversal to the leaves. Thus to make
a~hamiltonian description we need to extend $\fConn$. This is also known as
\emph{time dependent hamiltonians}.

\subsubsection{A symplectic extension of $\fConn$}
A \emph{level-$\divisor$ unipotent structure\/} on a~principal $G$-bundle $E$
is a~reduction of $E|_{\divisor}$ to $U$, where we view $\divisor$ as
a~closed subscheme of~$X$. Such a~structure~$\eta$ gives rise to a~Borel
structure (i.e.\ to a reduction of $E|_{\divisor}$ to $B$), which we denote
by $\eta_\fb$.

Let $\overline{\Conn}_U$ be the moduli space of triples $(E,\nabla,\eta)$,
where $(E,\nabla)\in\overline{\Conn}$, $\eta$ is a~unipotent level-$\divisor$
structure such that $\eta_\fb$ is \emph{compatible with $\nabla$}. In other
words, this structure is a~trivialization of $E$ at $x_i$ up to the order
$n_i-1$ (for all $i$) with the requirement that the coefficients of the polar
part of $\nabla$ are in~$\fb$. Two trivializations are considered the same if
they differ by an element of $U(O_\divisor)$.

Denote by $\Conn_U$ the open subspace of $\overline{\Conn}_U$ parameterizing
triples $(E,\nabla,\eta)$ with an additional condition that $\nabla$ is
non-resonant.

We shall construct a~natural map $\nu$ from $\Conn_U$ to $\fConn$. Take any
triple $(E,\nabla,\eta)\in\Conn_U$, let $x_i$ be an irregular point. Let
$\tilde\eta$ be any trivialization of $E$ at~$x_i$ up to the order $n_i-1$
such that $\tilde\eta$ extends $\eta$. Then $\tilde\eta$ gives rise to
a~framing~$\eta_0$ of $E$ at $x_i$. Let $A$ be the coefficient of $\nabla$ at
$z_i^{-n_i}dz_i$ relative to the framing $\eta_0$. Then $A\in\fb\cap\fgr$,
thus there is a~unique $u\in U$ such that $\Ad_uA\in\fh$. Then $u\eta_0$ is
a~unique framing at $x_i$ compatible with both $\nabla$ and~$\eta$, this
gives the desired map:
\[
    (E,\nabla,\eta)\mapsto(E,\nabla,u\eta_0).
\]

If follows from Proposition~\ref{modulistack} that $\overline{\Conn}_U$ is an
algebraic stack, while $\Conn_U$ is an algebraic space.

\begin{proposition}\label{connu}
$\Conn_U$ is a~smooth algebraic space with a natural symplectic structure.
\end{proposition}

This proposition will be proved in \S\ref{sehamappr}.

\subsubsection{Isostokes Hamiltonians}
Let $\conn_n$ be the scheme (of infinite type) of non-resonant connections
with a~pole of order $n$ on the trivial $G$-bundle on formal disk. We define
$\conn_n^B$ as the~subscheme of $\conn_n$, consisting of connections with
polar part in~$B$. We define $G_n^U$ as the group of loops of the form
\[
\exp(u_0+u_1z+\ldots+u_{n-1}z^{n-1}+g_nz^n+\ldots),
\]
where $u_i\in\fu$ for $i=0,\ldots,n-1$, $g_i\in\fg$ for $i\ge n$.

Take $(E,\nabla,\eta)\in\Conn_U$. Restricting $\nabla$ to the formal
neighborhood of~$x_i$, we get a~singular connection on the trivialized
punctured disk. Since $E$ is reduced to~$U$ up to the order $n_i-1$ at $x_i$
and $\nabla$ is compatible with this reduction, we get an element of
$\conn^B_{n_i}/G^U_{n_i}$. Thus we get a~map:
\[
IT_U:\Conn_U\to\prod_{i:n_i\ge2}\conn_{n_i}^B/G^U_{n_i}.
\]
The target of this map is a~Poisson variety, indeed, the $i$-th multiple is
an open subset in the hamiltonian reduction at 0 of the space of all
connections on the formal punctured disk with respect to $G^U_{n_i}$ (we
shall see in~\S\ref{hamdim} that the target of~$IT_U$ is a smooth affine
variety).

We come to the following commutative diagram:
\begin{equation}\label{CD}
\begin{CD}
\Conn_U @>IT_U>> \prod\limits_{i:n_i\ge2}\left(\conn^B_{n_i}/G^U_{n_i}\right)\\
@V\nu VV @VVV\\
\fConn @>IT>> (\fhr)^{l_{irr}}\times(\fh)^{\deg\divisor-l-l_{irr}}
\end{CD}
\end{equation}

We call a~tangent vector $v$ to $\Conn_U$ \emph{isostokes\/} if $\nu_*v$ is
isostokes (see Definition~\ref{isostokesvectors}).

\begin{theorem}\label{hamappr}
\setcounter{noindnum}{0}\noindstep\label{hamappr1} The map $IT_U$ is
a~Poisson map. If a~hamiltonian on $\Conn_U$ factors through~$IT_U$, then the
corresponding hamiltonian vector field on $\Conn_U$ is isostokes.\\
\noindstep\label{hamappr2} This construction gives the whole isostokes
deformation in the following sense: if~$v$ is a~tangent vector to the target
of $IT$ (see (\ref{CD})) at the point $IT(\nu(E,\nabla,\eta))$, where
$(E,\nabla,\eta)\in\Conn_U$, then there is a~hamiltonian vector field $v_H$,
whose hamiltonian factors through $IT_U$, such that $v_H(E,\nabla,\eta)$
projects to $v$.

\end{theorem}
We shall give an heuristic proof of the part (\ref{hamappr1})
in~\S\ref{sheur} and a~rigorous proof in~\S\ref{sehamappr}. The part
(\ref{hamappr2}) will be proved in~\S\ref{hamdim}. We shall also see
in~\S\ref{hamdim} that the number of linearly independent (at a~generic
point) hamiltonian vector fields produced by the part (\ref{hamappr1}) of the
theorem is equal to the dimension of the isostokes distribution on $\Conn_U$
given by~(\ref{dimfol}).

\begin{remark}
If the residues of $\nabla$ at regular singular points are in $\fgr$, we can
think about unipotent level-$\divisor$ structures as follows: reduce
$E|_\divisor$ to a~$B$-bundle. Under the condition that $\nabla$ preserves
this $B$-structure, this reduction is unique up to the action of $l$ copies
of the Weyl group. This $B$-bundle gives rise to a~$B/U$-bundle. The
unipotent level-$\divisor$ structure is a discrete choice of a $B$-reduction
plus a trivialization of the $B/U$-bundle over $\divisor$. Thus the dimension
of a~generic fiber of $\nu$ is
\begin{equation}\label{fnudim}
    (\deg\divisor-l_{irr})\rk\fg
\end{equation}
(recall that we do not have framings at regular points).

Another approach is to trivialize $E$ up to the order $n_i$ at every singular
point (this is called level-$\divisor$ structure). In that way we also obtain
a~smooth symplectic extension of $\Conn$. We decided on using unipotent
level-$\divisor$ structures because the dimension of a generic fiber of $\nu$
is equal to the codimension of the symplectic leaf of $\fConn$. Thus
$\Conn_U$ is a minimal symplectic extension of $\fConn$.
\end{remark}

\subsubsection{Dimensions}
The dimension of the analytic isostokes deformation is given by the dimension
of the target of $IT$, this is equal to
\begin{equation}\label{fdim}
    (\deg\divisor-l)\rk\fg.
\end{equation}

To calculate the dimension of the isostokes distribution on $\Conn_U$ we need
to add up (\ref{fdim}), (\ref{fnudim}), and $l_{irr}\rk\fg$, where the last
term comes from the deformations changing frames. The answer is
\begin{equation}\label{dimfol}
    (2\deg\divisor-l)\rk\fg.
\end{equation}

\section{Analytic isostokes deformations}\label{sanis}
In this section we shall give precise definitions of analytic isostokes
deformations and infinitesimal analytic isostokes deformations. Our primary
reference is~\cite{Boalch}.

\subsection{Stokes solutions and multipliers}\label{1.3}
Consider $(E,\nabla,s)\in\fConn$. Let $x_i$ be an irregular point, recall
that $z_i$ is an analytic coordinate at $x_i$. We can assume that
$z_i(x_i)=0$. Let $\mathcal U_i$ be the neighbourhood of $x_i$ given by
$|z_i|<\rho_i$ and $\mathcal V_i$ be given by $|z_i|<2\rho_i$ for some
$\rho_i>0$. The disks $\mathcal U_i$ and $\mathcal V_i$ will be fixed
throughout~\S\ref{sanis} and~\S\ref{sprfalgan}. We can assume that $\mathcal
V_i$ are disjoint.

For every irregular point $x_i$ we shall define the Stokes solutions and the
Stokes multipliers of $(E,\nabla,s)$ at $x_i$. Let us emphasize that the
discs $\mathcal U_i$ and $\mathcal V_i$ are defined only for
\emph{irregular\/} points $x_i$. We fix an irregular point $x_i$ for the
whole of~\S\ref{1.3}. Set $n=n_i$, $z=z_i$ for brevity.

\subsubsection{Stokes solutions}
Consider the coefficient $h=h_n\in\fhr$ at the leading term of the formal
normal form (\ref{fnf}) of $\nabla$. Let $\alpha$ be a~root of $\fg$ relative
to $\fh$. An \emph{anti-Stokes direction\/} corresponding to $\alpha$ is
a~ray in $\C$, emerging from the origin, on which $\alpha(h)z^{1-n}$ is real
and negative.

Let $r_1$ and $r_2$ be consecutive anti-Stokes directions. A \emph{Stokes
sector\/} is a~sector with the vertex at $x_i$ bounded by the directions
$r_1-\frac\pi{2n-2}$ and $r_2+\frac\pi{2n-2}$. We choose some Stokes sector
$S_1$ and then enumerate all the other Stokes sectors counterclockwise: $S_2,
S_3, \ldots, S_m$.

Note that the Stokes sectors cover $\mathcal V_i$ and that angular size of
each sector is greater than $\frac\pi{n-1}$. Notice also, that a single
anti-Stokes direction can correspond to more than one root, thus the number
of Stokes sectors can be different for irregular points of the same order.

We can write (\ref{fnf}) in the following form:
\begin{equation}\label{fnf2}
\nabla=d+dQ(z)+h_1\frac{dz}z,
\end{equation}
where $Q(z)$ is an $\fh$-valued polynomial in $\frac1z$. Then
$\Upsilon(z)=\exp(-h_1\ln z-Q(z))$ is a formal solution of $\nabla$, it is
called \emph{the canonical formal solution\/} of $\nabla$, $\Upsilon(z)$ will
be thought as a~multi-valued section of the trivial $G$-bundle over $\mathcal
V_i\setminus x_i$. We shall make it single-valued on every sector in the
following way: choose a~branch of $\Upsilon(z)$ over $S_1$ and subsequently
extend it to $S_2, S_3,\ldots, S_m$. Notice that these choices do not agree
on the intersection of $S_1$ and $S_m$.

Choose any trivialization of $E$ over $\mathcal V_i$, compatible with $s_i$.
Let $F(z)$ be a~unique formal gauge transformation such that $F(z)=1+O(z)$
and $F(z)$ takes $\nabla$ into the formal normal form~(\ref{fnf2}). Consider
the Stokes sector $S_j$. The \emph{Stokes solution\/} $\Phi_j(z)$ in $S_j$ of
$\nabla$ is a~unique solution of $\nabla$ in $S_j$ determined by the
requirement that the asymptotic expansion of $\Upsilon(z)^{-1}\Phi_j(z)$ at
the origin coincides with that of $F(z)$ (see~\cite{Boalch}, Theorem~1). This
solution of $\nabla$ does not depend on the trivialization chosen. However,
it depends on a choice of a~branch of $\Upsilon(z)$ and a choice of
a~numeration of Stokes sectors (i.e.\ a choice of a first Stokes sector).

\subsubsection{Stokes multipliers and the analytic classification of
connections} Let $S_j$ and $S_{j+1}$ be a~pair of consecutive Stokes sectors
($j\ne m$). The \emph{Stokes multiplier\/} for this pair of sectors is
$\Phi_j(z)\Phi_{j+1}(z)^{-1}\in G$. The Stokes multiplier, corresponding to
$S_m$ and $S_1$, is $e^{-2\pi ih_1}\Phi_m(z)\Phi_1(z)^{-1}\in G$. Notice that
Stokes multipliers are constant functions of $z$, since any two solutions of
a connection differ by a~constant element of $G$. It is well known (at least
for $G=GL$) that the formal normal form together with the Stokes multipliers
constitute a complete set of invariants of the local analytic classification
of connections, see Theorem~2 of~\cite{Boalch}. The Stokes multipliers belong
to certain subgroups of $G$ that are Weyl conjugate to $U$
(see~\cite{Boalch}, Lemma~6). This explains a somewhat peculiar formula for
the Stokes multiplier corresponding to~$S_m$ and $S_1$.

\subsubsection{Stokes solutions for families of connections}
Recall that Stokes solutions depend on a choice of a branch of $\Upsilon(z)$
and of a first Stokes sector. Let $(E(t),\nabla(t),s(t))$ be a non-resonant
family of connections over a~polydisk $\Delta\ni0$ (see Convention). Let us
make the choices above for $(E(0),\nabla(0), s(0))$. Let $\check S$ be
a~sector whose closure is contained in $S_j$ (it is assumed that the sectors
$S_j$ and $\check S$ have the same vertex). Then, shrinking $\Delta$ if
necessary and making a~continuous choice of a first Stokes sector for
$(E(t),\nabla(t),s(t))$, we can assume that the closure of~$\check S$ is in
the $j$-th Stokes sector for all $t\in\Delta$. We can also make a~continuous
choice of a~branch of $\Upsilon(z)$ for all $t\in\Delta$. Now we get a~family
of Stokes solutions $\Phi_j(z,t)$ defined on $\check S\times\Delta$. It
depends analytically on $z$ and $t$, see~Lemma~7 of~\cite{Boalch}. We shall
often omit $z$ in the notation below thus denoting this family by
$\Phi_j(t)$.

\subsection{Analytic isostokes deformations}

Let $(E(t),\nabla(t),s(t))$ be a non-resonant family. Denote the $j$-th
Stokes sector at $x_i$ by $S_j^i(t)$. We assume that $S_j^i(t)$ depends
continuously on $t$ (see above). Set $\CircX=X\setminus(\cup\mathcal U_i)$.

\begin{proposition}\label{analyticisostokes}
Given an analytic map $f$ from a~polydisk $\Delta$ to the target space of
$IT$ and $(E,\nabla,s)\in\fConn$, satisfying~(\ref{chp}), there is a~unique
up to isomorphism extension of $(E,\nabla,s)$ to a~non-resonant family
$(E(t),\nabla(t), s(t))$ over $\Delta$ such that\setcounter{noindnum}{0}\\
\noindstep\ $IT(E(t),\nabla(t),s(t))=f(t)$ for all $t\in\Delta$.\\
\noindstep\label{analiso2} The restriction of $(E(t),\nabla(t),s(t))$ onto
$\CircX$ is analytically constant.\\
\noindstep\label{analiso3} Let the closure of a~sector $\check S$ be
contained in $S^i_j(t)$ for all $t$. Let $\Delta'\subset\Delta$ be small
enough to define Stokes solution $\Phi^i_j(t)$ on $\check S\times\Delta'$,
then the restriction of $\Phi^i_j(t)$ to $(\partial\mathcal
U_i\cap\nobreak\check S)\times\Delta'$ does not depend on $t$ as a~section of
$E(t)|_{(\partial\mathcal U_i\cap\check S)\times\Delta'}$. (this condition
makes sense due to the condition~(\ref{analiso2})).
\end{proposition}

Notice that (\ref{analiso3}) is stronger than the~requirement that Stokes
multipliers do not change in the family.
\begin{proof}
To simplify notation we restrict to the case of a~single irregular
point~$x_1$. The general case is completely similar.

Denote by $\mathcal M$ the moduli space of triples $(\tilde
E,\tilde\nabla,\tilde s)$, where $\tilde E$ is a~principal $G$-bundle over
$\mathcal V_1$, $\tilde\nabla$ is a~non-resonant connection on $E$ with the
only pole at $x_1$ of the order $n_1$, $\tilde s$ is a~compatible framing at
$x_1$. Taking irregular type gives a~map $\mathcal
M\to\fhr\times(\fh)^{n_1-2}$. This is a~fibred bundle with a canonical flat
connection, obtained by deforming the irregular type, while preserving the
Stokes data. This is explained in~\cite{Boalch}, where this is called
\emph{the isomonodromic connection} (some details are given for $n_1=2$ only
but the general case is completely similar). We prefer to call this
connection \emph{the local isostokes connection}.

With this at hand we can finish the proof of the proposition. Let $(\tilde
E,\tilde\nabla,\tilde s)$ be the restriction of $(E,\nabla,s)$ to $\mathcal
V_1$. Then we use the local isostokes connection to extend $(\tilde
E,\tilde\nabla,\tilde s)$ to a~family $(\tilde E(t),\tilde\nabla(t),\tilde
s(t))$ of connections on $\mathcal V_1$ such that the irregular type of
$(\tilde E(t),\tilde\nabla(t),\tilde s(t))$ is $f(t)$. Let $(\hat
E,\hat\nabla)$ be the restriction of $(E,\nabla)$ to $X\setminus\mathcal
U_1$. It remains to patch $(\tilde E(t),\tilde\nabla(t))$ and $(\hat
E,\hat\nabla)$ together on $(\mathcal V_1\setminus\mathcal U_1)\times\Delta$.
The condition (\ref{analiso3}) of the proposition gives a~unique way to make
such a~patch.

In more detail, let $S^1_1(t),\ldots,S^1_m(t)$ be all the Stokes sectors for
$x_1$. Shrinking $\Delta$ if necessary we can choose a~system of sectors
$\check S_j$ such that (a) the closure of $\check S_j$ is contained in
$S^1_j(t)$ for all $t$ and (b) $\mathcal V_1\subset\cup\check S_j$. We have a
natural identification of $\hat E$ and $\tilde E$ over $(\mathcal
V_1\setminus\mathcal U_1)\times\{0\}$ and we use the condition
(\ref{analiso3}) to extend it to an identification of $\hat E$ and $\tilde E$
over $((\mathcal V_1\setminus\mathcal U_1)\cap\check S_j)\times\Delta$ for
every $j$. These identifications agree on the intersections, since the Stokes
multipliers do not change in $(\tilde E(t),\tilde\nabla(t),\tilde s(t))$. The
identifications respect the connection because $\Phi^1_j$ is a~solution of
$\nabla$.

It is clear that the way we patched $\hat E$ and $\tilde E$ is the only way
that satisfies~(\ref{analiso3}), thus the uniqueness. \end{proof}

\subsection{Infinitesimal isostokes deformations}\label{sinis}
Unfortunately, we do not know whe\-ther $\fConn$ is an algebraic scheme.
Therefore we shall use somewhat oblique way to define the algebraic isostokes
deformation. The problem is that $\fConn$ parameterizes algebraic families of
connections, while there are no algebraic isostokes families of connections
parameterized by smooth varieties. Thus we shall introduce the notion of
\emph{infinitesimally isostokes family\/} of connections.

Let $(E(t),\nabla(t),s(t))$ be an isostokes family of connections,
parameterized by a~smooth manifold $\Delta$. The restriction of $E(t)$ onto
$\CircX\times\Delta$ can be trivialized locally over $\Delta$. Indeed, if
$\Delta'\subset\Delta$ is an analytic disk, then $\CircX\times\Delta'$ is
a~Stein manifold and the claim follows from the Oka--Grauert principle,
see~\cite{Grauert}.

Fix such a~trivialization. Then the restriction of $\nabla(t)$ onto
$\CircX\times \Delta'$ becomes a~family of $\fg$-valued 1-forms on $\CircX$.
Thus the condition (\ref{analiso2}) of Proposition~\ref{analyticisostokes}
reads as follows: there is a~$\fg$-valued 1-form $\nabla$ on~$\CircX$ and
a~family of $G$-valued functions~$R(t)$ (where $t\in\Delta'$) on~$\CircX$
such that $\nabla(t)=\Ad_{R(t)}\nabla$. Here `$\Ad$' is the natural action of
$G$-valued functions on connections by gauge transformations. Below we shall
also use the infinitesimal action of $\fg$-valued functions on connections,
which we denote by $\ad$.

The restrictions of the Stokes solutions to $\partial\mathcal U_i$ can be
also viewed as $G$-valued functions in this trivialization. Then the
condition (\ref{analiso3}) of the definition becomes: the restriction of
$\Phi^i_j(t)R(t)$ to $(\partial\mathcal U_i\cap\check S)\times\Delta'$ does
not depend on~$t$.

As was mentioned in~\S\ref{ssalgap}, it will be convenient for us to work
with an extended version of the isostokes deformation. We want to add the
deformations that do not change connections but change framings at irregular
points. If we change a framing~$s_i$ to $Cs_i$, where $C\in G$, then the
Stokes solution $\Phi_j^i$ transforms into
$C^{-1}\Phi_j^i$. Thus we get a~weaker version of (\ref{analiso3}):\\

\vskip-3pt

(\ref{analiso3}$'$) For all $i$ there exists a~family $C_i(t)$ of elements of
$T$ such that the restriction of $C_i(t)^{-1}\Phi^i_j(t)R(t)$ to
$(\partial\mathcal U_i\cap\check S)\times\Delta'$, where $\check S$ as in the
part~(\ref{analiso3}) of Proposition~\ref{analyticisostokes}, does not depend
on $t$.

It is easy to write the infinitesimal version of these conditions.

\begin{definition}\label{infisostokes}
A non-resonant family $(E(t),\nabla(t),s(t))$ over a smooth manifold $\Delta$
is called \emph{infinitesimally isostokes at $t_0\in\Delta$ in the direction
of $v\in T_{t_0}\Delta$\/} if after trivializing the restriction of $E(t)$
onto $\CircX\times\Delta'$ (for some small neighbourhood $\Delta'$ of $t_0$)
we can find a~$\fg$-valued function $R$ on $\CircX$ and for each irregular
$x_i$ an element $c_i\in\fh$
such that\setcounter{noindnum}{0}\\
\noindstep\label{infis1} $L_v\nabla(t)=\ad_R\nabla(t_0)$ on $\CircX$\smcn\\
\noindstep\label{infis2} $L_v\Phi^i_j(t)=-\Phi^i_j(t_0)R+c_i\Phi^i_j(t_0)$ on
$\partial\mathcal U_i\cap S^i_j$ for all $i$ and $j$. Here $L_v$ is the
directional derivative in the direction of $v$ at $t_0$, $\Phi^i_j(t_0)R$ is
the usual left shift of $R$ on the tangent bundle of $G$, $c_i\Phi^i_j(t_0)$
is the right shift of $c_i$.
\end{definition}

It is clear that an isostokes family is infinitesimally isostokes at every
point and in every direction.

\section{Proof of Theorem~\ref{algan}}\label{sprfalgan}

Consider an algebraic non-resonant family $(E(t),\nabla(t),s(t))$ over
a~smooth variety $\Delta\ni t_0$. It gives rise to a~map $f:\Delta\to\fConn$.
Choose $p\in X$, $p\notin\cup\mathcal V_i$.

\subsection{Proof of the part `If' of the theorem}\label{iffirstpart}
Suppose that $v\in T_{t_0}\Delta$ is such that~$f_*v$ is an isostokes vector.
We need to show that $(E(t),\nabla(t),s(t))$ is an infinitesimally isostokes
family in the direction of $v$.

There is an \'etale neighbourhood $\iota:\Delta'\to\Delta$ of $t_0$ such that
$E(t)|_{\OpenX\times\Delta'}$ is trivial and $E(t)|_{(X\setminus p)\times
\Delta'}$ is trivial, since every $G$-bundle over a~family of affine curves
is trivial locally over the base in the \'etale topology,
see~\cite{Sorger,Drinfeld}. It is enough to show that the restriction of
$(E(t),\nabla(t),s(t))$ to $\Delta'$ is infinitesimally isostokes in the
direction of $(\iota_*)^{-1}v$, since every \'etale morphism is a~local
analytic diffeomorphism. Thus we can assume from the beginning that the
restrictions of $E(t)$ to $\OpenX\times\Delta$ and to $(X\setminus p)\times
\Delta$ are trivial.

Let us trivialize the restriction of $E(t)$ onto $\OpenX\times\Delta$. In
this trivialization the restriction of $\nabla(t)$ becomes a~family of
$\fg$-valued 1-forms, denote it by $\hat\nabla(t)$. Then we can re-write the
definition of $f_*v$ being isostokes in the following form: there is
a~$\fg$-valued (algebraic) function $R$ on $\OpenX$ such that
\begin{equation}\label{fisovect}
L_v\hat\nabla(t)=\ad_R\hat\nabla(t_0).
\end{equation}
We shall have to work in a neighbourhood of $\divisor$, thus we need
a~trivialization of~$E(t)$ in this neighbourhood. To this end we trivialize
$E(t)$ on $(X\setminus p)\times\Delta$. In this trivialization $\nabla(t)$ is
again a~$\fg$-valued 1-form, denote it by $\tilde\nabla(t)$. These two
trivializations (over $\OpenX\times\Delta$ and over $(X\setminus p)\times
\Delta$) are related by a~transition function $Z(t):(\OpenX\setminus p)\times
\Delta\to G$, $t\in\Delta$.

We identify $G$ with some subgroup of $GL$ via any exact representation, this
will simplify calculations. The condition (\ref{infis1}) in
Definition~\ref{infisostokes} is obvious: indeed, we just restrict $R$ in
(\ref{fisovect}) from $\OpenX$ to $\CircX$. Thus it suffices to check the
condition~(\ref{infis2}). In the trivialization of $E(t)$ on $(X\setminus
p)\times\Delta$ we get
\begin{equation}\label{fisovect2}
L_v\tilde\nabla(t)=\ad_P\tilde\nabla(t_0),
\end{equation}
where
\begin{equation}\label{fP}
P=(L_vZ(t))Z(t_0)^{-1}+Z(t_0)RZ(t_0)^{-1}.
\end{equation}

Fix an irregular point $x_i$. Recall that $z_i$ is an analytic coordinate on
$\mathcal V_i$ such that $z(x_i)=0$. We restrict $P$ and $\tilde\nabla(t)$
onto $\mathcal V_i$ and $\mathcal V_i\times\Delta$ respectively. From now on
we shall be working on $\mathcal V_i$, since the statement we need to prove
depends solely on the restrictions of our objects to this disc. We can also
assume that $t_0=0$. We shall work in the analytic setup, thus we view
$\Delta$ as an analytic manifold. Moreover, we can assume that $\Delta$ is
a~disk in~$\C$ (indeed, first we reduce to the case when $\Delta$ is a
polydisk, then we take the appropriate 1-dimensional section of this
polydisk). We emphasize that our objects depend on $z_i$, which we omit in
the notation. Finally, we set $z=z_i$, $n=n_i$, $\mathcal U=\mathcal U_i$,
and $\mathcal V=\mathcal V_i$ for brevity.

Write
\[
\tilde\nabla(t)=d+\sum_{j=-n}^\infty A_j(t)z^jdz.
\]

\begin{lemma}\label{assume}
Changing the trivialization of $E(t)$ on $\mathcal V\times\Delta$ by an
analytic gauge transformation we can assume that\setcounter{noindnum}{0}\\
\noindstep\label{assume1} this trivialization is compatible with the framing
$s_i$\smcn\\
\noindstep\label{assume2} $A_j(t)\in\fh$ for $-n\le j\le n$ and all $t$\smcn\\
\noindstep\label{assume3} $P$ is a polynomial in $\frac1z$\smcn\\
\noindstep all the coefficients of $P$ are in $\fh$.
\end{lemma}
\begin{proof}
Clearly, we can assume that (\ref{assume1}) is satisfied. Every connection
can be brought to its formal normal form up to any power of $z$ by an
analytic (even a polynomial) gauge action. It is easy to see that this can be
done in a family. Thus (\ref{assume2}) is clear. Note that since we have a
compatibly framed connection, this gauge change can be taken of the form
$1+O(z)$. It follows from~(\ref{fisovect2}) that the coefficients of $P$ are
in $\fh$ up to the coefficient at $z^{2n}$. Indeed, otherwise RHS
of~(\ref{fisovect2}) would not be in $\fh$ up to the order $n$, since
$\tilde\nabla(0)$ is non-resonant. Write
\[
P=P_-+P_+,
\]
where $P_-$ is a polynomial in $\frac1z$, $P_+$ is a polynomial in $z$
without the constant term. We can assume that $t$ is a coordinate on $\Delta$
such that $v=\frac{\partial}{\partial t}$. Change the trivialization of
$E(t)$ on $\mathcal V\times\Delta$ by means of $\exp(-tP_+)$. Then $P$
changes to (see~(\ref{fP})):
\[
P-P_+=P_-.
\]
Thus we get~(\ref{assume3}). The condition~(\ref{assume1}) of the lemma is
not corrupted by this trivialization change since $P_+$ has no constant term
and the condition~(\ref{assume2}) is not corrupted, since the coefficients of
$-tP_+$ are in $\fh$ up to the order $2n$. \end{proof}

\subsubsection{The infinitesimal change of the canonical formal solution} Write
\[
\tilde\nabla(t)=d+d_zQ(t)+\Lambda(t)\frac{dz}z+O(1),
\]
where $d_z$ is the differential with respect to $z$,
$Q(t)=\sum\limits_{j=1-n}^{-1}\frac{A_{j-1}(t)}jz^j$. By (\ref{assume2}) of
Lemma~\ref{assume}, the formal normal form of $\tilde\nabla(t)$ is just its
polar part, denote it by $\tilde\nabla_0(t)$.

Let
\[
\Upsilon(t)=\exp(-\Lambda(t)\ln z-Q(t))
\]
be the canonical formal solution of $\tilde\nabla(t)$. We want to study how
$\Upsilon(t)$ changes in the direction of $v$. Set
\[
\Theta=\Upsilon(0)^{-1}(L_v\Upsilon(t))=-L_v Q(t)
\]
(here we use that $L_v\Lambda(t)=0$, which follows from~(\ref{fisovect2})).
Let $P_0$ be the constant term of~$P$, we claim that
\begin{equation}\label{theta}
\Theta=P_0-P.
\end{equation}
Indeed, the definition of $\Theta$ shows that $-d_z\Theta$ is the polar part
of
\[
L_v\tilde\nabla(t)=\ad_P\tilde\nabla=d_zP+[P,A(0)],
\]
where we have written $\tilde\nabla(t)=d+A(t)dz$. This polar part is equal to
$d_zP$ by Lemma~\ref{assume}. Thus $-\Theta$ and $P$ differ by a term which
does not depend on $z$ but $\Theta$ is a polynomial in $\frac1z$ without the
constant term,  and the claim follows.

\subsubsection{The infinitesimal change of the Stokes solutions}
Recall the notion of Stokes solutions. The disk $\mathcal V$ is covered by
$m$ Stokes sectors. Let $F(t)=1+O(z)$ be the formal series in $z$ taking
$\tilde\nabla(t)$ into its formal normal form. According to
Lemma~\ref{assume}, we can assume that $F(t)=1+O(z^{n+1})$. Fix a Stokes
sector $S$ for $\tilde\nabla(0)$. Take a sector $\check S$ of angular size
greater than $\frac\pi{n-1}$ whose closure is in $S$. As before, we may
assume that the corresponding Stokes solution $\Phi(t)$ is defined on $\check
S\times\Delta$. This solution of~$\tilde\nabla(t)$ is uniquely determined by
the requirement that the asymptotic expansion of $\Upsilon(t)^{-1}\Phi(t)$ in
$z$ at the origin coincides with one of $F(t)$ for all $t$.

\begin{lemma}\label{lhalf}
\begin{equation}\label{fphiprime}
L_v\Phi=-\Phi(0)P+c\Phi(0),
\end{equation}
where $c\in\fh$ is a constant matrix.
\end{lemma}
\begin{proof}
We shall prove this by a direct computation. Set $\Psi=L_v\Phi$. We have
\[
d_z\Phi(t)=-\Phi(t)A(t)dz.
\]
Applying $L_v$ to both sides of this equation, we get a variation equation
for~$\Psi$:
\[
d_z\Psi=-\Psi A(0)dz+\Phi(0)(d_zP+[P,A(0)]\,dz).
\]
It is easy to verify that $-\Phi(0)P$ satisfies the same differential
equation. Thus $\Psi+\Phi(0)P$ is a solution of $\tilde\nabla(0)$, which
gives
\begin{equation}\label{apriori}
\Psi=-\Phi(0)P+c\Phi(0),
\end{equation}
where $c$ is a constant matrix. It remains to show that $c\in\fh$. We have
\begin{equation}\label{asympt}
\Upsilon(t)^{-1}\Phi(t)\sim F(t)
\end{equation}
for every $t$. Moreover, it follows from the proof of Lemma~7
of~\cite{Boalch} that this asymptotic expansion is uniform in some
neighbourhood of the origin. Thus it follows from the Cauchy integral formula
that we can apply $L_v$ to both sides of this equation: we get
\begin{equation}\label{f32}
\Upsilon(0)^{-1}\Psi-\Theta\Upsilon(0)^{-1}\Phi(0)\sim L_vF(t).
\end{equation}
Substituting $\Psi$ from~(\ref{apriori}) and using~(\ref{theta})
and~(\ref{asympt}), we get
\[
-F(0)P+\Upsilon(0)^{-1}c\Phi(0)+PF(0)-P_0F(0)\sim L_vF(t).
\]
Since $F(0)=1+O(z^{n+1})$ and $L_vF(t)=O(z^{n+1})$, we get that
\[
\Upsilon(0)^{-1}c\Upsilon(0)=P_0+O(z).
\]
Now, let $\mathop{\oplus}\limits_\alpha\fg_\alpha\oplus\fh$ be the root
decomposition of $\fg$. Suppose that the projection of $c$ to $\fg_\alpha$ is
not zero for some $\alpha$. Denoting the corresponding character of~$T$ by
$\exp(\alpha)$, we see that $\exp(\alpha)(\Upsilon(0))$ must be bounded in
$\check S$. However, this function is of the form $z^\lambda e^{f(z)}$, where
$f$ is a polynomial in $\frac1z$ of degree $n-1$. Such a function cannot be
bounded in $\check S$, since the angular size of $\check S$ is greater than
$\frac\pi{n-1}$. This contradiction shows that $c\in\fh$. \end{proof}

In the trivialization of $E(t)$ over $\OpenX\times\Delta$ the Stokes solution
$\Phi(t)$ transforms into $\hat\Phi(t)=Z(t)^{-1}\Phi(t)$, and
(\ref{fphiprime}) becomes
\[
    L_v\hat\Phi(t)=-\hat\Phi(0)R+c\hat\Phi(0),
\]
and we see that $(E(t),\nabla(t),s(t))$ is infinitesimally isostokes in the
direction of $v$.

\subsection{Proof of the part `Only if' of the theorem}\label{4.2}
This is, in some sense, a rearrangement of the previous proof. Let us assume
that $(E(t),\nabla(t),s(t))$ is infinitesimally isostokes in the direction of
$v$. Again, we can, after passing to an \'etale neighbourhood of
$t_0\in\Delta$ assume that the restrictions of $E(t)$ to $\OpenX\times
\Delta$ and $(X\setminus p)\times\Delta$ are trivial.

Then we have to show that there is an \emph{algebraic\/} $\fg$-valued
function~$R$ on~$\OpenX$ such that~(\ref{fisovect}) holds on $\OpenX$. The
condition (\ref{infis1}) of Definition~\ref{infisostokes} gives an
\emph{analytic\/} function~$R$ such that~(\ref{fisovect}) holds on $\CircX$.
Shrinking $\Delta$ if necessary, we can choose a system of sectors~$\check
S^i_j$ such that (a) the closure of $\check S^i_j$ is contained in $S^i_j(t)$
for all $t$ and (b)~$\mathcal V_i$ is contained in $\cup_j\check S^i_j$.

We extend $R$ to $\OpenX$ solving the equation in the
Definition~\ref{infisostokes}, i.e.\ setting
\[
    R=\Phi^i_j(0)^{-1}(c\Phi^i_j(0)-L_v\Phi^i_j(t))
\]
on $\mathcal V_i\cap\check S^i_j$. The condition~(\ref{infis2}) of
Definition~\ref{infisostokes} together with the analytic continuation
principle show that these definitions of~$R$ agree on $\mathcal
V_i\cap\CircX\cap\check S^i_j$. The compatibility on $\check S^i_j\cap\check
S^i_{j+1}$ follows from the compatibility on $\mathcal
V_i\cap\CircX\cap\check S^i_j\cap\check S^i_{j+1}$ by analytic continuation.
It remains to show that $R$ does not have essential singularities.

This is easy to check, employing parts of the proof above. First, it is
enough to check that $P$ given by~(\ref{fP}) has no essential singularities,
since~$Z(t)$ is meromorphic. In $\mathcal U_i$ we have~(\ref{fphiprime}) with
some $c\in\fg$. We also have~(\ref{f32}). Substituting the former into the
latter we see that
\[
    P\sim F(0)^{-1}cF(0)-F(0)^{-1}\Theta F(0)-F(0)^{-1}L_vF(t)
\]
and therefore $P$ has no essential singularity at $x_i$. Theorem~\ref{algan}
is proved.

\section{Algebraic structures on moduli spaces of connections}\label{secstack}
In this section we shall prove Proposition~\ref{modulistack}.

Let $\mathrm{Bun}_GX$ be the moduli stack of principal $G$-bundles on $X$.
This is an algebraic stack, locally of finite type over $\C$,
see~\cite{Sorger}, Corollary~3.6.6, see also~\cite{Laumon},~4.14.2.1 for the
case $G=GL$. Thus we can restrict ourselves to families over locally
noetherian schemes below.

Let $\F$ be a~divisor on $X$, denote by $\overline{\Conn}_\F$ the lax functor
(or 2-functor) from affine schemes to groupoids defined by:
\[
    S\mapsto\{\mbox{pairs $(E(t),\nabla(t))$}\}+\{isomorphisms\}.
\]
Here $E(t)$ is a~$G$-bundle on $X\times S$, $\nabla(t)$ is a~connection along
$X$ with the pole divisor bounded by $\F\times S$. This is a~stack in the
\'etale topology, the proof is essentially the same as in the case of
$\mathrm{Bun}_GX$ (see~\cite{Deligne}, Theorem~4.5). The only additional
ingredient is that connections can be glued in the \'etale topology, this is
obvious.

Now we shall prove that $\overline{\Conn}=\overline{\Conn}_\divisor$ is an
algebraic stack. Choose an ample divisor $\E$ on $X$. Let
$\mathrm{Bun}^{(k)}_GX$ be the stack of $G$-bundles $E$ such that
\[
    H^1(X,\ad E\otimes\Omega^1(\divisor+k\E))=0.
\]
Precisely, this stack parameterizes bundles $E(t)$ over $X\times S$ ($t\in
S$) with
\begin{equation}\label{vanishing}
    R^1p_*(\ad E(t)\otimes\Omega((\divisor+k\E)\times S))=0,
\end{equation}
where $p:X\times S\to S$ is the natural projection, $\Omega=\Omega^1_{X\times
S/S}$ is the sheaf of relative differentials. Clearly, this is an open (and
hence algebraic) substack in $\mathrm{Bun}_GX$, it is also clear that
$\mathrm{Bun}_GX=\cup_k \mathrm{Bun}^{(k)}_GX$ (see~\cite{Laumon},~4.14.2.1
for more details).

Let $\Conn_{\divisor+k\E}^{(k)}$ be the substack of
$\overline{\Conn}_{\divisor+k\E}$ parameterizing pairs $(E,\nabla)$ such that
$E\in\mathrm{Bun}^{(k)}_GX$. Consider the forgetful 1-morphism of stacks
\[
    \lambda:\Conn_{\divisor+k\E}^{(k)}\to\mathrm{Bun}^{(k)}_GX.
\]
We claim that it is representable. Indeed, let $S\to\mathrm{Bun}^{(k)}_GX$ be
any morphism. It corresponds to a~$G$-bundle $E(t)$ over $X\times S$. This
bundle has a~connection with poles on $(\divisor+k\E)\times S$ \'etale
locally over $S$, since the local obstruction to the existence of such
a~connection is in the vanishing sheaf~(\ref{vanishing}).

Thus (locally over $S$) the set of connections on $E(t)$ is identified with
the total space of
\[
    R^0p_*(\ad E(t)\times\Omega((\divisor+k\E)\times S)).
\]
By~(\ref{vanishing}) and the Riemann--Roch theorem, this sheaf is locally
free. It follows that the fiber
\[
    S\times_{\mathrm{Bun}^{(k)}_GX}\Conn_{\divisor+k\E}^{(k)}
\]
is an affine bundle over $S$, hence, a scheme.

Thus $\lambda$ is representable, therefore $\Conn_{\divisor+k\E}^{(k)}$ is an
algebraic stack. It follows that the substack of $\overline{\Conn}_\divisor$
corresponding to the bundles satisfying~(\ref{vanishing}) is algebraic as
well, since it is a~closed substack of $\Conn_{\divisor+k\E}^{(k)}$. Thus
$\overline{\Conn}$ is a~union of an increasing sequence of open algebraic
substacks, hence algebraic.

It remains to show that $\overfConn$ is an algebraic space. The stack of
connections with arbitrary (not necessarily compatible) framings is algebraic
since the forgetful 1-morphism to $\overline{\Conn}$ is representable. Thus
$\overfConn$ is an algebraic stack, since it is a~closed substack of the
latter stack. However, framed connections do not possess automorphisms, thus
$\overfConn$ is an algebraic space.

This completes the proof of Proposition~\ref{modulistack}.

\section{The double quotient construction and isostokes hamiltonians}\label{sheur}

In this section we give an heuristic proof of the part~(\ref{hamappr1}) of
Theorem~\ref{hamappr}. It is based on an infinite-dimensional symplectic
reduction. Unfortunately, there are some technical difficulties in such an
approach, therefore we give another proof in the next section. The current
proof explains how the theorem was invented making clear the connection with
the paper~\cite{FrenkelBenZvi}.

\subsection{The double quotient construction}\label{doubleq}
Let $G((z))$ be the group of $G$-valued functions on the punctured formal
disk (the loop group). Denote by $LG$ the group $\prod_1^l G((z))$. We can
identify $LG$ with the group of $G$-valued functions on the formal
neighbourhood of $\divisor$. Let
\[
L_+G=\prod_{i=1}^lG^U_{n_i}\subset LG
\]
be the subgroup of ``positive loops''. Let $L_XG=G(\OpenX)$ be the group of
$G$-valued functions on $\OpenX$. Such a function can be restricted to the
formal neighbourhood of $\divisor$, which gives an embedding
$L_XG\hookrightarrow LG$. Then the stack of $G$-bundles with level-$\divisor$
unipotent structures is isomorphic to the double quotient
\[
L_XG\setminus LG/L_+G.
\]
The similar statement is well known for the stack of $G$-bundles without
additional structures, see~\cite[Theorem~5.1.1]{Sorger} and
\cite[Theorem~4.1.1]{FrenkelBenZvi}. In our case the proof is completely
similar. Morally, an element of $LG$ is viewed as a $G$-bundle, trivialized
over both $\OpenX$ and the formal neighbourhood of $\divisor$, while the
factoring by $L_XG$ and~$L_+G$ amounts to forgetting these trivializations.

Notice that $LG$ and $L_XG$ are ind-groups, $L_+G$ is an affine group of
infinite type, see~\cite{Sorger}, \S8.

\subsection{Generalities on hamiltonian quotients}\label{ssgen}
Notice that we use the expressions ``hamiltonian quotient'' and ``hamiltonian
reduction'' as synonyms.

Let $Y$ be a Poisson ind-scheme, equipped with a hamiltonian action of an
ind-group $K$. Let $\mathbf{O}$ be a coadjoint orbit in $\Lie(K)^*$. Denote
by $Y//_\mathbf{O}K$ the hamiltonian quotient of~$Y$ by $K$ at $\mathbf{O}$.
This is the quotient $\mu^{-1}(\mathbf{O})/K$, where $\mu$ is the moment map.

Let $H:Y//_\mathbf{O}K\to\C$ be a function. We can lift it to a function
\[
    \tilde H:\mu^{-1}(\mathbf{O})\to\C.
\]
By a \emph{lift of $H$ to $Y$\/} we mean a function $\hat H:Y\to\C$ such that
its restriction to $\mu^{-1}(\mathbf{O})$ coincides with $\tilde H$. Clearly,
such a lift is not unique. Actually, we can start with a function $H$,
defined on an open subset of $Y//_\mathbf{O}K$, then its lift is a function
on an open subset of $Y$. Note that in our case ``open'' means open in
\`etale topology, this is where the difficulties come from.

Denote by $v_H$ the hamiltonian vector field corresponding to $H$, by
$\{\cdot,\cdot\}$ the Poisson bracket.

\begin{lemma}\label{sympred}
Let $H$ and $H_1$ be functions on an open subset of $Y//_\mathbf{O}K$, $\hat
H$ and~$\hat H_1$
be their lifts to $Y$. Then\setcounter{noindnum}{0}\\
\noindstep $\{\hat H,\hat H_1\}$ is a lift of $\{H,H_1\}$.\\
\noindstep A vector field $v_{\hat H}$ is tangent to $\mu^{-1}(\mathbf{O})$,
the restriction of $v_{\hat H}$ to $\mu^{-1}(\mathbf{O})$ is
$K$-equi\-va\-ri\-ant and descends to $v_H$.
\end{lemma}
These are standard hamiltonian reduction facts.

\subsection{Double quotient presentation of $\overline{\Conn}_U$}\label{isodiscuss}

Let $T^*_1LG$ be a twisted cotangent bundle to $LG$, parameterizing pairs
$(g,\nabla)$, where $g\in LG$ and $\nabla$ is a connection on the formal
punctured neighbourhood of~$\divisor$. Let
$\conn=\lim\limits_\rightarrow\conn_n$ (see~\S\ref{shamappr}) be the
ind-scheme of connections on the trivial formal punctured disk. We may view
$T_1^*LG$ as the space parameterizing $G$-bundles on $X$, trivialized both
on~$\OpenX$ and the formal neighbourhood of $\divisor$ with a singular
connection $\nabla$ on this formal neighbourhood (compare
with~\S\ref{doubleq}). We may view $\nabla$ as an element of $(\conn)^l$,
using any of two trivializations of $T^*_1LG$, which gives two isomorphisms
\[
    T^*_1LG\cong LG\times(\conn)^l.
\]
Denote the corresponding projections $T^*_1LG\to(\conn)^l$ by $p^R$ and
$p^L$. Precisely,~$p^R$ corresponds to the trivialization of the bundle over
the open curve $\OpenX$, while $p^L$ corresponds to the trivialization over
the formal neighbourhood of $\divisor$. Clearly, the adjoint action of $LG$
on $(\conn)^l$ intertwines these projections.

This twist can be also explained by symplectic reduction. Let $\hat\fg$ be an
affine Kac-Moody algebra that is the canonical central extension of the loop
algebra $\fg((z))$. Then $\tilde\fg=\prod_{i=1}^l\hat\fg$ is a central
extension of $\Lie(LG)$. Let $\tilde G$ be the corresponding central
extension of $LG$, then $\Lie\tilde G=\tilde\fg$. The center of $\tilde\fg$
integrates to a central subgroup of $\tilde G$, isomorphic to
$(\C^\times)^l$. This gives a $(\C^\times)^l$ action on $T^*\tilde G$. We
have
\begin{equation}\label{fredat1}
    T_1^*LG=T^*\tilde G//_{\mathbf 1}(\C^\times)^l,
\end{equation}
where $//_{\mathbf 1}$ states for the symplectic reduction at
${\mathbf1}=(1,\ldots,1)\in(\Lie(\C^\times)^l)^*$.

The proof of Theorem~\ref{hamappr} is based on the following presentation:
\begin{equation}\label{fsympdq}
    \overline{\Conn}_U=L_XG\RightHamRed T_1^*LG//_0L_+G.
\end{equation}
This formula also gives the desired symplectic structure on $\Conn_U$. Notice
that this double quotient can be thought as a single symplectic quotient with
respect to the group $L_XG\times L_+G$.

We call $v\in T_1^*LG$ an \emph{isostokes\/} vector, if $p^R_*v=0$. This
definition, clearly, agrees with~(\ref{fsympdq}) and
Definition~\ref{isostokesvectors} in the following sense: suppose $v$ is a
vector, tangent to the zero-level of the moment map and such that its
projection to $\overline{\Conn}_U$ is tangent to $\Conn_U$. Then it is
isostokes iff its projection to $\Conn_U$ is.

\begin{remark}
We have isomorphisms of ind-schemes:
\[
    T^*_1LG\cong T^*LG \cong LG\times\Lie(LG)^*.
\]
However, each of these spaces is equipped with both left and right actions
of~$G$. The isomorphisms can be chosen either left or right equivariant but
not bi-equivariant. Thus these spaces are different as $G$-bimodules. In
addition, first two spaces carry different symplectic structures.
\end{remark}

\subsection{Completion of the ``proof'' of the first part of
Theorem~\ref{hamappr}}\label{Sscomp} The target of $IT_U$ can be identified
with an open subset of
\[
(\conn)^{l_{irr}}\Bigl/\Bigl/_0\Bigl(
\textstyle\prod\limits_{\:i:n_i\ge2}G^U_{n_i}\Bigr).
\]
Let $H$ be a hamiltonian on $\Conn_U$ that factors through $IT_U$: $H=f\circ
IT_U$. Let~$\bar f$ be a lift of $f$ to $(\conn)^{l_{irr}}$ (recall, that
this is a function on an open subset of $(\conn)^{l_{irr}}$). Using the
natural projection $(\conn)^l\to(\conn)^{l_{irr}}$ we can lift $\bar f$ to a
function $\hat f$ on an open subset of $(\conn)^l$. Then it is easy to see
that $\hat H=\hat f\circ p^L$ is a lift of $H$ with respect
to~(\ref{fsympdq}).

Taking into account Lemma~\ref{sympred} and the discussion in the end of
\S\ref{isodiscuss}, we see that our theorem reduces to the following
statement:\\
\vskip-6pt
\centerline{\emph{$p^L$ is a Poisson map, $p^R_*v_{\hat H}=0$.}}

\vskip3pt
This is not hard to verify by a direct calculation. However, we can reduce it
to some general nonsense, using~(\ref{fredat1}). Indeed, let $\tilde p^L$ and
$\tilde p^R$ be the projections $T^*\tilde G\to\tilde\fg^*$, corresponding to
the left and right trivializations of the cotangent bundle.

We can identify $(\conn)^l$ with a subspace in $\tilde\fg^*$ as usual. Now,
let us lift $\hat f$ to $\tilde f:\tilde\fg^*\to\C$ and set $\tilde H=\tilde
f\circ\tilde p^L$. Applying Lemma~\ref{sympred} again we reduce the theorem
to the following
statement:\\
\centerline{\emph{$\tilde p^L$ is a Poisson map, $\tilde p^R_*v_{\tilde
H}=0$.}}

\vskip3pt
This is true for any Lie group $K$: the left projection $T^*K\to\Lie(K)^*$ is
Poisson; hamiltonians that factor through this projection preserve the leaves
of the right trivialization of $T^*K$.


\def\U{\mathcal U}
\def\Supp{\mathop{\mathrm{Supp}}}
\def\higgs{\mathop{\mathrm{Higgs}}}
\def\K{\mathcal K^\bullet}
\def\Ker{\mathop{\mathrm{Ker}}}
\def\Image{\mathop{\mathrm{Image}}}
\def\Spf{\mathop{\mathrm{Spf}}}

\section{The symplectic structure on $\Conn_U$ via hypercohomology}\label{sehamappr}

In this section we give a rigorous proof of the part (\ref{hamappr1}) of
Theorem~\ref{hamappr}.

\subsection{Tangent Space to $\overline{\Conn}_U$.}

The following presentation of the tangent space to the stack
$\overline{\Conn}$ is well known:
\[
T_{(E,\nabla)}\overline{\Conn}=\Hyper^1(X,\ad E\xrightarrow{\ad_\nabla}\ad
E\otimes\omega(-\divisor)).
\]
Here $\omega=\Omega^1(X)$ is the canonical bundle on $X$. We are going to use
not the formula above but its version for the tangent space to
$\overline{\Conn}_U$ at $(E,\nabla,\eta)$. Denote by $\ad(E,\eta)\subset\ad
E$ the sheaf of infinitesimal automorphisms of $E$ preserving $\eta$.  Its
stock at $x\notin\Supp\divisor$ coincides with the one of $\ad E$, while its
stock at $x_i\in\Supp\divisor$ is (non-canonically) isomorphic to the set of
loops of the form
\[
    g_0+g_1z+g_2z^2+\ldots,
\]
where $g_j\in\fu$ for $0\le j<n_i$. Denote by $\higgs(E,\eta)$ the sheaf of
$(\ad E)$-valued 1-forms with polar part bounded by $\divisor$ that are
compatible with $\eta$.

\begin{proposition}
There is a canonical isomorphism:
\begin{equation}\label{ftanspace}
T_{(E,\nabla,\eta)}\overline{\Conn}_U\cong
\Hyper^1(X,\ad(E,\eta)\xrightarrow{\ad_\nabla}
\higgs(E,\eta)).
\end{equation}
\end{proposition}

\begin{proof}
Consider an affine open cover $\U_\alpha$ of $X$ such that
\[
\U_\alpha\cap\U_\beta\cap\Supp(\divisor)=\emptyset
\]
for all $\alpha\ne\beta$. Consider a section $s_\alpha$ of $E$ over each of
$\U_\alpha$ such that $s_\alpha$ agrees with $\eta$. The transition functions
between $s_\alpha$ form a $G$-valued 1-cocycle $\phi_{\alpha\beta}$. In the
trivialization $s_\alpha$ the connection $\nabla$ becomes a $\fg$-valued
1-form on $\U_\alpha$, denote it by~$\theta_\alpha$. The pair
$(\phi_{\alpha\beta},\theta_\alpha)$ determines the triple $(E,\nabla,\eta)$
up to an isomorphism.

Denote the complex in~(\ref{ftanspace}) by $\K$ and consider its Czech
resolution with respect to~$\mathcal U_\alpha$:
\[
\begin{CD}
@>>> 0 @>>>\scriptstyle\ad(E,\eta) @>\ad_\nabla>>
\scriptstyle\higgs(E,\eta) @>>>0 @>>>\\
@. @VVV @VVV @VVV @VVV\\
@>>> 0 @>>>\scriptstyle C^0(\ad(E,\eta))
@>>>\genfrac{}{}{0pt}{}{C^1(\ad(E,\eta))\oplus}{C^0(\higgs(E,\eta))} @>>>
\genfrac{}{}{0pt}{}{C^2(\ad(E,\eta))\oplus}{C^1(\higgs(E,\eta))} @>>>
\end{CD}
\]
The complex at the bottom is the cone of the morphism
\[
\ad_\nabla:C^\bullet(\ad(E,\eta))\longrightarrow C^\bullet(\higgs(E,\eta)).
\]

Suppose that we have an infinitesimal deformation
\[
\phi_{\alpha\beta}\mapsto\phi_{\alpha\beta}\exp(\epsilon\psi_{\alpha\beta}),\qquad
\theta_\alpha\mapsto\theta_\alpha+\epsilon\nu_\alpha.
\]
Changing the trivializations $s_\alpha$ we can check by a direct computation
that the pair $(\psi_{\alpha\beta},\nu_\alpha)$ is naturally identified with
an element of
\[
C^1(\ad(E,\eta))\oplus C^0(\higgs(E,\eta)).
\]
It is also easy to see that the compatibility condition on intersections is
equivalent to $(\psi_{\alpha\beta},\nu_\alpha)$ being a cocycle of $\K$. To
conclude the proof of the proposition it remains to show that two cocycles
give isomorphic deformations iff they differ by a~coboundary. This can be
also done by a direct computation. \end{proof}

\subsection{Smoothness of $\Conn_U$}
It is a standard fact, that the obstruction to smoothness of $\Conn_U$ is in
$\Hyper^2(X,\K)$. Note that $\K$ is a self-dual complex, thus by the
Grothendieck duality $\Hyper^2(X,\K)$ is dual to $\Hyper^0(X,\K)$. The latter
space vanishes. Indeed, framed connections have no automorphisms and there is
an algebraic map $\nu:\Conn_U\to\fConn$. Thus $\Conn_U$ is a smooth algebraic
space.

\subsection{The symplectic structure on $\Conn_U$}
Since $\K$ is a self-dual complex, the Grothendieck duality gives a
non-degenerate 2-form~$\varpi$ on $\Conn_U$. We need to check that this
2-form is closed. Let us write the explicit formulae first. Suppose that we
have two tangent vectors represented by cocycles
\[
(\psi^i_{\alpha\beta},\nu^i_\alpha)\in C^1(\ad(E,\eta))\oplus
C^0(\higgs(E,\eta)),\qquad i=1,2.
\]
Then the value of the symplectic form on these vectors is
\[
    \Res\langle\psi^1_{\alpha\beta}\nu^2_\alpha-
    \psi^2_{\alpha\beta}\nu^1_\alpha\rangle,
\]
where $\langle\cdot,\cdot\rangle$ is the Cartan-Killing form,
$\Res:H^1(X,\omega)\to\C$ is the natural isomorphism.

Let $v_x, v_y, v_z$ be in $T_{(E,\nabla,\eta)}\Conn_U$, we need to check that
$d\varpi(v_x,v_y,v_z)=0$. Since $\Conn_U$ is smooth, we can find a map
$F:\Spf[[x,y,z]]\to\Conn_U$ such that $F_*\frac d{dx}=v_x$, $F_*\frac
d{dy}=v_y$, $F_*\frac d{dz}=v_z$ ($\Spf$ states for formal spectrum). It is
enough to show that
\[
    d(F^*\varpi)\textstyle\bigl(\frac d{dx},\frac d{dy},\frac d{dz}\bigr)=0
\]
at the unique closed point of $\Spf[[x,y,z]]$. Let us calculate $F^*\varpi$
in the coordinates $x$, $y$, $z$. The map $F$ corresponds to a family of
connections parameterized by $\Spf[[x,y,z]]$. We can write it by a cocycle as
before:
\begin{equation*}
\begin{split}
\mathcal
C=(&\psi^{\alpha\beta}+\psi^{\alpha\beta}_xx+\psi^{\alpha\beta}_yy+\psi^{\alpha\beta}_zz+
\psi^{\alpha\beta}_{xy}xy+\psi^{\alpha\beta}_{yz}yz+\psi^{\alpha\beta}_{xz}xz+...,\\
&\nu^{\alpha}+\nu^{\alpha}_xx+\nu^{\alpha}_yy+\nu^{\alpha}_zz+
\nu^{\alpha}_{xy}xy+\nu^{\alpha}_{yz}yz+\nu^{\alpha}_{xz}xz+...)
\end{split}
\end{equation*}
(we omit irrelevant terms). Consider the vector field $\frac d{dx}$ on
$\Spf\C[[x,y,z]]$; it gives rise to a vector field on $\Conn_U$ along
$\Spf\C[[x,y,z]]$, i.e.\ a map $\Spf\C[[x,y,z]]\times\dual\to\Conn_U$ given
by the following cocycle:
\[
 \mathcal C+\epsilon(\psi_x^{\alpha\beta}+
 \psi_{xy}^{\alpha\beta}y+\psi_{xz}^{\alpha\beta}z+\ldots,\;\nu_x^{\alpha}+
 \nu_{xy}^{\alpha}y+\nu_{xz}^{\alpha}z+\ldots).
\]
We have similar vector fields for $\frac d{dy}$ and $\frac d{dz}$. Thus, up
to irrelevant terms we have
\[
\scriptstyle F^*\varpi\bigl(\textstyle\frac d{dx},\frac d{dy}\bigr)=
\scriptstyle\varpi((\psi_x^{\alpha\beta}+
 \psi_{xy}^{\alpha\beta}y+\psi_{xz}^{\alpha\beta}z,\nu_x^{\alpha}+
 \nu_{xy}^{\alpha}y+\nu_{xz}^{\alpha}z),(\psi_y^{\alpha\beta}+
 \psi_{xy}^{\alpha\beta}x+\psi_{yz}^{\alpha\beta}z,
 \nu_y^{\alpha}+\nu_{xy}^{\alpha}x+\nu_{yz}^{\alpha}z)).
\]
We see that
\begin{multline*}
\textstyle\left.\frac d{dz}\right|_{x=y=z=0}F^*\varpi\bigl(\frac d{dx},\frac
d{dy}\bigr)=\\
\Res(\langle\psi_x^{\alpha\beta},\nu_{yz}^\alpha\rangle-
\langle\psi_{yz}^{\alpha\beta},\nu_x^{\alpha}\rangle+
\langle\psi_{xz}^{\alpha\beta},\nu_y^\alpha\rangle-
\langle\psi_y^{\alpha\beta},\nu_{xz}^{\alpha}\rangle).
\end{multline*}
It remains to write out two similar expressions obtained from this by a
cyclic permutation of $x$, $y$ and $z$ and add them up. They add up to zero,
this concludes the proof of Proposition~\ref{connu}.

\subsection{Isostokes vectors}
Using the description of the tangent space above, it is easy to see that
$v\in T_{(E,\nabla,\eta)}\Conn_U$ is an isostokes vector iff it is in the
kernel of the natural map
\[
\Hyper^1(X,\K)\to\Hyper^1(\OpenX,\K)=\Hyper^1(X,j_*j^*\K),
\]
where $j:\OpenX\hookrightarrow X$ is the natural embedding. One may think
about this as about the tangent map to the restriction map from connections
on~$X$ to connections on~$\OpenX$.

Clearly,
\[
    j_*j^*\K=\lim_{N\to+\infty}\K\otimes O(N\divisor).
\]
Thus the space of isostokes vectors can be written as
\begin{equation}\label{fisohyper}
    \lim_{N\to+\infty}\,\Ker(\Hyper^1(X,\K)\to\Hyper^1(X,\K\otimes O(N\divisor))).
\end{equation}

\subsection{Hamiltonians that factor through $IT_U$}
Let $H:\Conn_U\to\C$ be a function that factors through $IT_U$. Then
$dH|_{(E,\nabla,\eta)}$ vanishes on the corresponding vertical subspace
$T_0\subset T_{(E,\nabla,\eta)}\Conn_U$. We want to describe this subspace
$T_0$ in terms of hypercohomology. Let \
\[
\divisor'=\sum_{i:n_i\ge2}n_ix_i
\]
be the irregular part of $\divisor$. We claim that
\begin{equation}\label{fvertical}
    T_0=\lim_{N\to+\infty}\Image(\Hyper^1(X,\K\otimes
    O(-N\divisor'))\to\Hyper^1(X,\K)).
\end{equation}

Indeed, suppose that $(\psi_{\alpha\beta},\nu_\alpha)\in T_0$. Look at the
Laurant expansion of $\theta_\alpha+\epsilon\nu_\alpha$ at an irregular point
$x_i\in\U_\alpha$. This expansion can be viewed as a tangent vector to
$\conn^B_{n_i}$ and this vector is tangent to an orbit of $G^U_{n_i}$ by the
definition of $T_0$.

It follows that for all $N$ there is an infinitesimal gauge transformation
$Z$ such that (1) $Z$ is defined on $\U_\alpha$, (2) $Z$ preserves $\eta$ and
(3) $\ad_Z\nu_\alpha$ vanishes up to the order $N$ at $x_i$. In other words,
we can make $\nu_\alpha$ vanish up to any order at irregular points by adding
a coboundary to $(\psi_{\alpha\beta},\nu_\alpha)$. And the claim follows.

\subsection{End of the proof of Theorem~\ref{hamappr}} It remains to prove that
the space~(\ref{fvertical}) contains the symplectic complement to the
space~(\ref{fisohyper}). Since $\divisor\succ\divisor'$, it is enough to show
that
\[
\Ker(\Hyper^1(X,\K)\to\Hyper^1(X,\K\otimes O(N\divisor)))
\]
is the symplectic complement to
\[
\Image(\Hyper^1(X,\K\otimes O(-N\divisor))\to\Hyper^1(X,\K))
\]
for all $N>0$. Let $i:\K\hookrightarrow\K\otimes O(N\divisor)$ be the natural
inclusion. Our statement follows from the functoriality of the Grothendieck
duality and the following commutative diagram
\[
\begin{CD}
(\K\otimes O(N\divisor))^{op}\otimes\omega @>i^{op}\otimes\mathrm{Id}>>
(\K)^{op}\otimes\omega\\
@VV\cong V @VV\cong V\\
\K\otimes O(-N\divisor) @>>>\K
\end{CD}
\]
The bottom map is the natural inclusion. This concludes the proof of the
first part of the theorem.

\section{Dimensions and an explicit construction of hamiltonians}\label{hamdim}

In this section we shall complete the proof of Theorem~\ref{hamappr}.

\begin{proposition}\label{uniqueform}
Every element of $\conn_n^B/G^U_n$ has a unique representative in $\conn_n^B$
of the form
\begin{equation}\label{NormalForm}
d+\alpha_n\frac{dz}{z^n}+\ldots+
\alpha_1\frac{dz}z+\beta_0dz+\ldots+\beta_{n-2}z^{n-2}dz,
\end{equation}
where $\alpha_i,\beta_i\in\fh$ for all $i$.
\end{proposition}

\begin{proof}
Take an element of  $\conn_n^B/G^U_n$, let
\[
d+\alpha_n\frac{dz}{z^n}+\ldots+\alpha_1\frac{dz}z+\beta_0dz+\beta_1zdz+\ldots
\]
be any of its representatives. We already have $\alpha_i\in\fb$ for all $i$.
We can put~$\alpha_n$ into~$\fh$ by the gauge action of a constant loop $g\in
U$.

Further, $[\alpha_n,\fu]=\fu$. It allows to put all the polar part of the
connection into $\fh$ by the gauge action of an appropriate loop of the form
\[
    \exp(u_1z+u_2z^2+\ldots+u_{n-1}z^{n-1}),
\]
where $u_i\in\fu$ for all $i$. Since $[\alpha_n,\fg]+\fh=\fg$, we can put the
terms of positive order into~$\fh$, using the gauge action of a loop of the
form
\[
    \exp(g_nz^n+g_{n+1}z^{n+1}+\ldots).
\]
To kill the terms of order higher than $n-2$, we use the appropriate loop of
the form
\[
    \exp(h_nz^n+h_{n+1}z^{n+1}+\ldots),
\]
where all $h$'s are in $\fh$. This proves the existence part of the
proposition. We leave the uniqueness to the reader. \end{proof}

\begin{corollary}
\setcounter{noindnum}{0}\noindstep\ For $n>1$ we have an isomorphism of
varieties:
\[
\conn_n^B/G^U_n\approx\fhr\times(\fh)^{2n-2}.
\]
\noindstep
\[
\dim\conn_n^B/G^U_n=(2n-1)\rk\fg.
\]
\end{corollary}

\begin{remark}
We see that the target space of~$IT_U$ is affine. Thus there are a lot of
global isostokes hamiltonians. These hamiltonians do not commute, since this
space has a non-trivial Poisson structure. It would be interesting to
construct commuting hamiltonians.
\end{remark}

\begin{proposition}\label{symple}
The symplectic leaves of $\conn_n^B/G^U_n$ are given by $\alpha_1=\const$,
see~(\ref{NormalForm}).
\end{proposition}
\begin{proof}
Let $A$ be an element of $\conn_n^B/G^U_n$ whose representative $\tilde A$ is
given by~(\ref{NormalForm}) (we use Proposition~\ref{uniqueform}). We shall
calculate the tangent space to the symplectic leaf
containing~(\ref{NormalForm}) in $\conn$. The Poisson structure on $\conn$
comes from the immersion $\conn\hookrightarrow\hat\fg^*$. Thus the tangent
space to the symplectic leaf at $\tilde A$ is given by $\ad_{\hat\fg}\tilde
A$. It is easy to see that it consists of exactly those $v\in T_{\tilde
A}\conn$ whose residue has no diagonal part (notice that $T_{\tilde
A}\conn=\fg((z))^*$).

Now, unwinding the definition of the hamiltonian reduction we obtain the
required statement. \end{proof}

Now we can prove the part (\ref{hamappr2}) of Theorem~\ref{hamappr}. Consider
the map
\[
    \phi_n:\conn_n^B/G^U_n\to\fhr\times(\fh)^{n-2}
\]
that assigns the $(n-1)$-tuple $(\alpha_n,\ldots,\alpha_2)$
to~(\ref{NormalForm}). Set
\[
\phi=\prod_{i:n_i\ge2}\phi_{n_i}
\]
(this is the right vertical arrow in~(\ref{CD})). We see that for every
tangent vector $v$ to the target of $IT$ at $\phi(A)$ there is a hamiltonian
$f$ on the target space of $IT_U$ such that $\phi_*(v_f|_A)=v$, where $v_f$
is the hamiltonian vector field corresponding to $f$. Taking $H=f\circ IT_U$
one completes the proof of Theorem~\ref{hamappr}.

\subsection{`Stupid' hamiltonians} Some of isostokes hamiltonians,
produced by Theorem~\ref{hamappr}, are `stupid': they do not change irregular
types of connections but the unipotent structures only (see the
diagram~(\ref{CD})). Here we shall describe these hamiltonians. According to
Proposition~\ref{uniqueform}, we can consider $\alpha_i$ for $i=1,\ldots,n$
and $\beta_i$ for $i=0,\ldots,n-2$ as the coordinates on $\conn_n^B/G^U_n$.
Define a coordinate $\alpha_i^j$ on the target of $IT_U$ as the composition
of $\alpha_i$ and the projection to the $j$-th multiple. Similarly we
define~$\beta_i^j$.

\begin{proposition}
$f\circ IT_U$ is a stupid hamiltonian iff $f$ does not depend on $\beta_i^j$.
\end{proposition}
\begin{proof}
Clearly, it suffices to prove that a hamiltonian $f:\conn_n^B/G^U_n\to\C$
satisfies $\phi_*v_f=0$ iff it does not depend on $\beta_i$'s.

Notice first, that
\[
    \{\alpha_i,\alpha_j\}=0
\]
for all $i$ and $j$, this easily follows from the presentation of
$\conn_n^B/G^U_n$ as a hamiltonian reduction of $\conn$ (compare with the
proof of Proposition~\ref{symple}). Thus $\phi$ is a~lagrangian fibration and
the claim follows. \end{proof}


\end{document}